\documentclass[letterpaper,11pt]{article}

 \hoffset= - 0.75 in

\title{Small Contingency Tables with Large Gaps}

\author{Seth Sullivant
 \\ {\small Department of Mathematics, University of California, Berkeley}}

\date{}

\usepackage{latexsym,array,delarray,amsthm,amssymb}
\usepackage[sumlimits]{amsmath}
\usepackage[dvips]{graphicx}

\theoremstyle{plain}
\newtheorem{thm}{Theorem}

\theoremstyle{definition}
\newtheorem{defn}[thm]{Definition}

\theoremstyle{remark}

\begin{document}
\maketitle
\begin{abstract}
We construct examples of contingency tables on $n$ binary random
variables where the gap between the linear programming lower/upper
bound and the true integer lower/upper bounds on cell entries is
exponentially large.   These examples provide evidence
that linear programming may not be an effective heuristic for
detecting disclosures when releasing margins of multi-way tables.
\end{abstract}

\section{Introduction}

A fundamental problem in data security is to determine what
information about individual survey respondents can be
inferred from the release of partial data.  The particular
instance of this problem we are interested in concerns the release
of margins of a multidimensional contingency table.  In
particular, given a collection of margins of a multi-way table,
can individual cell entries in the table be inferred.  This type
of problem arises when statistical agencies like a census bureau
release summary data to the public, but are required by law to
maintain the privacy of individual respondents.

Many authors \cite{B,Ch,C} have proposed that an individual cell
entry is secure if, among all contingency tables with the given
fixed marginal totals, the upper bound and lower bound for the
cell entry are far enough apart.  In general, solving the integer
program associated with finding the sharp integer upper and lower
bounds a cell entry is known to be NP-hard.  A heuristic which
has been suggested for approximating these upper and lower bounds
is to solve the appropriate linear programming relaxation.  Based
on theoretical results for 2-way tables and practical experience
for some small multi-way tables, some authors have suggested that
the linear programming bounds and other heuristics should always
constitute good approximations to the true bounds for cell values.

In this paper, we attempt to refute the claim that the linear
programming bounds are, in general, good approximations to the
true integer bounds.  In particular, we will show the following:

\begin{thm}
There is a sequence of hierarchical models on $n$ binary random
variable and a
collection of margins such that the gap between the linear
programming lower (upper) bounds and the integer programming lower
(upper) bounds for a cell entry grows exponentially in $n$.
\end{thm}

For instance, on 10 binary random variables, our construction
produces an instance where this difference is more than 100.  This
constitutes a significant discrepancy between the heuristic and
reality, in a problem of size which is quite small from the
practical standpoint.

The outline of this paper is as follows.  In the next section we
review hierarchical models and the algebraic techniques that we will
use to construct our examples.  The third section is devoted to
the explicit construction, and in the fourth section we discuss
practical consequences of our examples.

\section{Graphical Models, Gr\"obner Bases, and Graver Bases}

A hierarchical model is given by a collection of subsets
$\Delta$ of the $n$-element set $[n] := \{1,2,\ldots, n\}$
together with an integer vector $d = (d_1, \ldots, d_n)$.  Without
loss of generality, we can take $\Delta$ to be a simplicial complex.  In the
setting of probabilistic inference, a hierarchical model is intended
to encode interactions between a collection of $n$ discrete random
variables:  the number of states is the $i$-th random variable is
$d_i$ and there is an interaction factor between the set of random
variables indexed by each $F \in \Delta$ (see, for example,
\cite{L} for an introduction).  From the standpoint of data
security, $n$ is the number of dimensions of a multi-way
contingency table, the $d_i$ represent the number of levels in
each dimension, and the elements $F \in \Delta$ are the particular
margins that are released.  For the rest of this paper $d_i = 2$
for all $i$; that is, we are considering \emph{dichotomous} tables
or \emph{binary} random variables.

Computing the $\Delta$-margins of a multi-way table
is a linear transformation.  We denote by $A_\Delta$ the matrix in
the standard basis that computes these margins.  Finding the minimum
value for a cell entry given the $\Delta$-margins $\mathbf{b}$ amounts
to solving
the following integer program, which we denote $IP_{\Delta}$:

$$ \min u_{\mathbf{0}} \mbox{ subject to } $$
$$ A_\Delta \mathbf{u} = \mathbf{b}, \mathbf{u} \geq 0, \mathbf{u}
\mbox{ integral}. $$

\noindent  The linear programming relaxation drops the
integrality condition.  We denote it by $LP_{\Delta}$:

$$ \min u_{\mathbf{0}} \mbox{ subject to } $$
$$ A_\Delta \mathbf{u} = \mathbf{b}, \mathbf{u} \geq 0. $$

\noindent  The integer programming gap  $gap_-(\Delta)$ is the
largest difference  between the optimal solution of $IP_\Delta$
and $LP_\Delta$ over all feasible marginals $\mathbf{b}$
\cite{HS}.  Explicitly computing the integer programming gap is a
difficult problem, even for quite small models $\Delta$.  However,
using properties of Gr\"obner bases, it is easy to give lower
bounds on this gap.  Recall the definition of a Gr\"obner basis:

\begin{defn}
A reduced Gr\"obner basis $G_\mathbf{c}$ of $A_\Delta$ with respect to
the cost vector 
$\mathbf{c}$ is a minimal set of improving vectors
that solves the integer 
program $IP_{\Delta,\mathbf{c}}$ for any feasible right hand side
$\mathbf{b}$.
\end{defn}

In the literature of discrete optimization, Gr\"obner bases are often
called test sets.
A lower bound on $gap_-(\Delta)$ is given by inspecting the
coordinates of the Gr\"obner basis with respect to the cost
vector $\mathbf{c} = \mathbf{e}_{00 \cdots 0}$.

\begin{thm}[\cite{HS}, Corollary 4.3]\label{thm:HS}
The value $gap_{-}(\Delta)$ is greater than or equal to one less than
the largest coordinate $g_{00 \cdots 0}$ of any element in the reduced
Gr\"obner basis $G_\mathbf{c}$ of $A_\Delta$.
\end{thm}

\noindent  The precise definition of the Gr\"obner bases can be found
in \cite{St}, however, we will restrict to a special family of models
where the Gr\"obner basis elements we need have a simpler description.
For this, we
will need to recall the definition of the Graver basis.  Note that any
integer vector $\mathbf{u}$, can be written uniquely as $\mathbf{u} =
\mathbf{u}^+ - \mathbf{u}^-$, where $\mathbf{u}^+$ and $\mathbf{u}^-$
are nonnegative with disjoint support.

\begin{defn}
A nonzero integer vector $\mathbf{u} \in \ker(A_\Delta)$ is called primitive is there
does not exist an integer vector $\mathbf{v} \in
\ker(A_\Delta) \setminus \{\mathbf{0}, \mathbf{u} \}$ such that
$\mathbf{v}^+ \leq \mathbf{u}^+$ and $\mathbf{v}^- \leq \mathbf{u}^-$.
The set of vectors $\{\mathbf{u} \in  \ker(A_\Delta) | \mathbf{u}
\mbox{ is primitive} \} $ is called the \emph{Graver basis} of
$A_\Delta$.
\end{defn}

Given a simplicial complex $\Gamma$ on $[n-1]$ there is a natural
construction of a new simplicial complex $\Delta = logit(\Gamma)$ on
$[n]$ which corresponds to taking the logit model with a binary
response variable.  The new model is defined as

$$logit(\Gamma) := \{ S \cup \{n\} | S \in \Gamma \} \cup 2^{[n-1]}$$

\noindent where $2^{[n-1]}$ is the set of all subsets of $[n-1]$.  Note
that $\ker(A_\Gamma)$ and $\ker(A_{logit(\Gamma)})$ are isomorphic,
and there is a natural identification:  $\mathbf{u} \in
\ker(A_\Gamma)$ if and only if $(\mathbf{u},-\mathbf{u}) \in \ker(
A_{logit(\Gamma)})$.  This follows by inspecting the condition
required by the margin associated to
the facet $[n-1]$ of $logit(\Gamma)$.  A fundamental fact about logit
models is that their Gr\"obner bases are easy to describe in terms of
the Graver basis of $A_\Gamma$, namely:

\begin{thm}[\cite{St} Theorem 7.1]\label{thm:gr}
Let $\Gamma$ be a model and $\Delta = logit(\Gamma)$ then:

\begin{enumerate}
\item $Gr(A_\Delta) = \{(\mathbf{u},-\mathbf{u})| \mathbf{u} \in
  Gr(A_\Gamma)\}$,
\item $\{ g \in Gr(A_\Delta)|\mathbf{c} \cdot \mathbf{g} > 0 \}
  \subseteq G_\mathbf{c}$. 
\end{enumerate}
\end{thm}

Note that Theorem \ref{thm:gr} is only true when the response
variable is binary.  We now have all the tools in hand to construct
our example.

\section{The Construction}

Our main result is the following:

\begin{thm}
For each $n \geq 3$, there is a hierarchical model $\Delta_n$ on $n$-binary
random variables such that $$gap_-(\Delta_n) \geq 2^{n-3} - 1.$$
\end{thm}

A similar statement about exponential growth of the gap for upper
bounds can be derived by an analogous arument.

\begin{proof}
Our strategy will be to construct
a hierarchical model $\Delta_n$ which has Gr\"obner basis elements
whose $\mathbf{0}$ entry is large.  This will force the large gap by
Theorem \ref{thm:HS}.

Let $\Gamma_n$ be the hierarchical model on $n-1$ random variables

$$\Gamma_n = \{ S | S \subset [n-2], S \neq [n-2] \} \cup \{ \{n-1 \} \}.$$

\noindent  That is, $\Gamma_n$ is the union of the boundary of an $n-3$
simplex together with an isolated point. Take $\Delta_n =
\mathrm{logit}(\Gamma_n)$.  To show the theorem with respect to
$\Delta_n$ is suffices to show that $A_{\Gamma_n}$ has elements in its
Graver basis that have large entries in their $\mathbf{0}$ coordinate,
by Theorem \ref{thm:gr}.

Consider the vector

$$\mathbf{f}_n = 2^{n-3}e_{(\mathbf{0},0)} + \sum_{\mathbf{i} | \mathbf{i}
  \neq \mathbf{0} , \sum i_j even } e_{(\mathbf{i}, 1)} \, \, - \, \,
  (2^{n-3}-1)e_{(\mathbf{0},1)} - \sum_{\mathbf{i}|\sum i_j odd}
  e_{(\mathbf{i},0)}.$$

Here $e_{(\mathbf{i},k)}$ denotes the standard unit vector whose index
is $(\mathbf{i},k) \in \{0,1\}^{n-1}$;  that is, $e_{(\mathbf{i},k)}$
is the integral table whose only nonzero entry is a one in the
$(\mathbf{i},k)$ position.  Note that $\mathbf{i} \in \{0,1\}^{n-2}$
is an index on the first $n-2$ random variables.

We will now show that $\mathbf{f}_n$ is a primitive vector in
$\ker(A_{\Gamma_n})$.  First we must show that $\mathbf{f}_n \in
\ker(A_{\Gamma_n})$; that is, the positive and the negative part of
$\mathbf{f}_n$ have the same margins with respect to $\Gamma_n$.  However, the
margin with respect to any of the subsets $S \subset [n-2], S \neq
[n-2]$ is the same: namely, it is the vector $\mathbf{m}_n$ given by
$$\mathbf{m}_n = (2^{n-3}-1)e_{\mathbf{0}} +  \sum_{\mathbf{i} \in
  \{0,1\}^{n-3}} e_\mathbf{i}. $$
\noindent  The margin with respect to $\{n-1\}$ is the vector
$\mathbf{m}'_n$ given by
$$\mathbf{m}'_n = 2^{n-3}e_0 + (2^{n-3}-1)e_1.$$
\noindent  In particular, these margins are the same and so
$\mathbf{f}_n$ belongs to $\ker(A_{\Gamma_n})$.

Now we must show that $\mathbf{f}_n$ is a primitive vector in
$\ker(A_{\Gamma_n})$.  Suppose to the contrary that there was some
nontrivial $\mathbf{g}_n \in \ker(A_{\Gamma_n})$  such that
$\mathbf{g}_n^+ \leq
\mathbf{f}_n^+$ and $\mathbf{g}_n^- \leq
\mathbf{f}_n^-$.  Suppose that one of the coordinates of
$\mathbf{g}_n^+$ was nonzero in a position indexed by some
$(\mathbf{i},1)$ with $\sum i_j$ even.  Then this forces
$\mathbf{g}_n^+$  to
have nonzero entries in all the possible positions indexed by $(\mathbf{i},1)$
with $\sum i_j$ even if the margins with respect to the $S \subset
[n-2]$ are to be the same in $\mathbf{g}_n^+$ and $\mathbf{g}_n^-$.
However, this implies that the margin of $\mathbf{g}_n^+$ with respect
to $\{n-1\}$ has an entry of $2^{n-3} -1$ in the $1$ position.  This
forces $\mathbf{g}_n = \mathbf{f}_n$ if $\mathbf{g}_n \in 
\ker(A_{\Gamma_n})$.  On the other hand, since $\mathbf{g}_n \neq \mathbf{0}$,
it must have some positive entry.  However, its only positive entry
could not be in the $(\mathbf{0},0)$ position since this would force
a negative entry in some position $(\mathbf{i},0)$.  By the preceding
argument, this implies that $\mathbf{g}_n = \mathbf{f}_n$ and thus
$\mathbf{f}_n$ is a primitive vector.

\end{proof}

To explicitly construct an example of a set of margins $\mathbf{b}$
with respect to $\Delta_n$ where the gap between the LP and IP optima
is $2^{n-3} - 1$ just take

$$\mathbf{u} =(2^{n-3}-1)e_{(\mathbf{0},0,0)} + \sum_{\mathbf{i} | \mathbf{i}
  \neq \mathbf{0} , \sum i_j even } e_{(\mathbf{i}, 1,0)} \, \, + \, \,
  (2^{n-3}-1)e_{(\mathbf{0},1,1)} + \sum_{\mathbf{i}|\sum i_j odd}
  e_{(\mathbf{i},0,1)},$$

\noindent and $\mathbf{b} = A_{\Delta_n} \mathbf{u}$.  It follows that
$\mathbf{u}$ cannot be improved to an nonnegative integer table with
smaller $(\mathbf{0},0,0)$ coordinate by appealing to the Gr\"obner
basis.  However, the nonnegative rational vector

$$\mathbf{v} = \mathbf{u} - \frac{2^{n-3}-1}{2^{n-3}} (\mathbf{f}_n,
-\mathbf{f}_n)$$

\noindent has the same margins $\mathbf{b}$ as $\mathbf{u}$ but has
$(\mathbf{0},0,0)$ coordinate $0$.

\section{Discussion}

In this paper, we constructed an example to show that the gap between
the linear programming lower bounds and the integer programming lower
bounds for a cell entry can be exponentially large in the number of
binary random variables of a hierarchical model.  Previous explicit
constructions of this type \cite{DS} gave gaps that were linear in the
number of random variables.  

There are a number of possible modifications to our result which can
be made, to produce examples of different flavors.  For instance,
small modifications of our argument can 
be used to produce exponential gaps between the linear programming and
integer programming upper bounds for cell entries.  Furthermore, by
adding extra dimensions by subdividing $\Delta$,  and using some of
the techniques in \cite{DS}, one can produces instances of purely
graphical models with these exponential growth properties.

While it is not clear how often, given a random collection of margins
$\mathbf{b}$, one should expect to encounter the exponentially large
gaps we have demonstrated, we expect that for problems on large sparse
tables, large gaps between the LP and IP solutions will be not be
exceptional.  This feeling is based on the observation that if any gap
value can occur, then so can all the integer values smaller than this
gap.  This suggests that research needs to be done to
determine better heuristics for approximating bounds on cell entries
in large sparse tables.

\end{document}